%% file: jour-13-acp-siam.tex
\documentclass[final,leqno]{siamltex}

\usepackage{amsmath,amssymb}
\usepackage{graphicx}
\usepackage{epstopdf}
\usepackage{subfigure}
\usepackage{caption}
\usepackage{setspace}
\usepackage{textcomp}
\usepackage{algorithmic}
\usepackage{algorithm}
\usepackage{cases}
\usepackage{bbm}
\usepackage{upgreek}
\usepackage[square,comma,numbers,sort&compress]{natbib}
\usepackage{hyperref}
\usepackage[all]{hypcap}

\hypersetup{colorlinks=true,linkcolor=blue,citecolor=red}

\input{def_siam}

\begin{document}

\setcounter{footnote}{1}
\title{A Convergence Proof of the Split Bregman Method for Regularized Least-Squares Problems}
\author{Hung Nien\!\!\!\!\and Jeffrey A. Fessler\thanks{Department of Electrical Engineering and Computer Science, University of Michigan, Ann Arbor, MI 48109, USA. ({\tt \{hungnien, fessler\}@umich.edu})}}

\maketitle

\begin{abstract}
The split Bregman (SB) method \cite{goldstein:09:tsb} is a fast splitting-based algorithm that solves image reconstruction problems with general $\ell_1$, e.g., total-variation (TV) and compressed sensing (CS), regularizations by introducing a single variable split to decouple the data-fitting term and the regularization term, yielding simple subproblems that are separable (or partially separable) and easy to minimize. Several convergence proofs have been proposed \cite{esser:09:aol,setzer:09:sba,cai:09:sbm}, and these proofs either impose a ``full column rank'' assumption to the split or assume exact updates in all subproblems. However, these assumptions are impractical in many applications such as parallel magnetic resonance (MR) and X-ray computed tomography (CT) image reconstructions \cite{goldstein:09:tsb,ramani:12:asb,chun:13:ecs,cauley:13:hsh}, where the inner least-squares problem usually cannot be solved efficiently due to the highly shift-variant Hessian. In this paper, we show that when the data-fitting term is quadratic, e.g., in image restoration problems with Gaussian noise, the SB method is a convergent alternating direction method of multipliers (ADMM) \cite{glowinski:75:slp,gabay:76:ada,eckstein:92:otd,afonso:11:aal}, and a straightforward convergence proof with inexact updates is given using \cite[Theorem 8]{eckstein:92:otd}. Furthermore, since the SB method is just a special case of an ADMM algorithm, it seems likely that the ADMM algorithm will be faster than the SB method if the augmented Largangian (AL) penalty parameters are selected appropriately. To have a concrete example, we conduct a convergence rate analysis of the ADMM algorithm with two split variables (the SB method is just a special case of the two-split ADMM algorithm) for image restoration problems with quadratic data-fitting term and regularization term. According to our analysis, we can show that the two-split ADMM algorithm can be faster than the SB method if the AL penalty parameter of the SB method is suboptimal. Numerical experiments were conducted to verify our analysis.
\end{abstract}

\section{Introduction} \label{sec:jour-13-acp:intro}
Consider a regularized least-squares optimization problem with a general convex regularizer:
\begin{equation} \label{eq:reprot-13-asc:ls_general_l1}
	\hat{\mb{x}}
	=
	\argmin
	{\mb{x}}
	{\ts\frac{1}{2}\norm{\mb{y}-\mb{Ax}}{2}^2+\fx{\Phi}{\msb{\Theta}\mb{x}}} \, ,
\end{equation}
where $\mb{y}$ is the noisy measurement, $\mb{A}$ is the system matrix, $\Phi$ is a convex potential function, and $\msb{\Theta}$ is an analysis matrix. For example, when $\Phi$ is the $\ell_1$-norm and $\msb{\Theta}$ is the discrete framelet transform matrix \cite{daubechies:03:fmb}, the regularized least-squares problem \eqref{eq:reprot-13-asc:ls_general_l1} is a frame-based image restoration problem \cite{cai:09:sbm}; when $\Phi$ is a smooth ``$\ell_1$-like'' potential function (such as the Huber function \cite{huber:81:rs,nikolova:05:aoh} and the Fair function \cite{fair:74:otr,fessler:99:cgp}) and $\msb{\Theta}$ is the finite difference matrix, the regularized least-squares problem \eqref{eq:reprot-13-asc:ls_general_l1} is an image restoration problem with an edge-preserving regularizer. To solve \eqref{eq:reprot-13-asc:ls_general_l1}, one can use the split Bregman (SB) method proposed by Goldstein \textit{et al.} \cite{goldstein:09:tsb}, which solves an equivalent constrained minimization problem:
\begin{equation} \label{eq:jour-13-acp:eq_ls_general_l1}
	\left(\hat{\mb{x}},\hat{\mb{v}}\right)
	=
	\argmin
	{\mb{x},\mb{v}}
	{\ts\frac{1}{2}\norm{\mb{y}-\mb{Ax}}{2}^2+\fx{\Phi}{\mb{v}}}
	\text{ s.t. }
	\mb{v}=\msb{\Theta}\mb{x}
\end{equation}
using the (alternating direction) augmented Lagrangian (AL) method. The iterates of the SB method are as follows:
\begin{equation} \label{eq:jour-13-acp:sb_iterates}
	\begin{cases}
	\iter{\mb{x}}{k+1}
	=
	\argmin
	{\mb{x}}
	{\ts\frac{1}{2}\norm{\mb{y}-\mb{Ax}}{2}^2+\ts\frac{\eta}{2}\norm{\msb{\Theta}\mb{x}-\iter{\mb{v}}{k}-\iter{\mb{e}}{k}}{2}^2} \\
	\iter{\mb{v}}{k+1}
	=
	\argmin
	{\mb{v}}
	{\fx{\Phi}{\mb{v}}+\ts\frac{\eta}{2}\norm{\msb{\Theta}\iter{\mb{x}}{k+1}-\mb{v}-\iter{\mb{e}}{k}}{2}^2} \\
	\iter{\mb{e}}{k+1}
	=
	\iter{\mb{e}}{k}-\msb{\Theta}\iter{\mb{x}}{k+1}+\iter{\mb{v}}{k+1} \, ,
	\end{cases}
\end{equation}
where the $\mb{x}$-update is a least-squares problem, and the $\mb{v}$-update is a proximal mapping of $\Phi$, which often can be solved efficiently, e.g., by soft-thresholding for the $\ell_1$ potential.

To prove the convergence of the SB method, Esser \cite{esser:09:aol} showed that the SB method is equivalent to the alternating direction method of multipliers (ADMM) \cite{glowinski:75:slp,gabay:76:ada,eckstein:92:otd,afonso:11:aal}, and Setzer \cite{setzer:09:sba} showed that the SB method can be interpeted as the Douglas-Rachford splitting (DRS) method \cite{douglas:56:otn,eckstein:92:otd,combettes:07:adr} applied to the dual problem. However, both \cite{esser:09:aol} and \cite{setzer:09:sba} assume that $\msb{\Theta}$ has full column rank, i.e., $\msb{\Theta}'\msb{\Theta}$ is \emph{invertible}, and show convergence proofs using \cite[Theorem 8]{eckstein:92:otd}. The full column rank condition holds when $\msb{\Theta}$ is a \emph{tight frame} as in frame-based image restoration problems. When $\msb{\Theta}$ is the finite difference matrix as in edge-preserving image restoration problems (and also in \cite{goldstein:09:tsb}), this assumption does not hold anymore, and the the proofs in \cite{esser:09:aol} and \cite{setzer:09:sba} are inapplicable. Differently, in \cite{cai:09:sbm}, assuming all the inner updates in \eqref{eq:jour-13-acp:sb_iterates} are exact, Cai \textit{et al.} proved the convergence of the SB method without using \cite[Theorem 8]{eckstein:92:otd} and therefore did not impose the ``full column rank'' assumption. In other words, the SB method is a convergent algorithm for any $\msb{\Theta}$ if \emph{all} the inner minimization problems in \eqref{eq:jour-13-acp:sb_iterates} are solved \emph{exactly}! Unfortunately, when some of the inner updates are inexact, e.g., the $\mb{x}$-update in parallel magnetic resonance (MR) and X-ray computed tomography (CT) image reconstructions \cite{goldstein:09:tsb,ramani:12:asb,chun:13:ecs,cauley:13:hsh}, we still lack convergence proofs of the SB method. In this paper, we first show the equivalence of the SB method and a convergent ADMM algorithm for solving \eqref{eq:reprot-13-asc:ls_general_l1}, and then give a simple convergence proof of the SB method that allows inexact updates when the data-fitting term is quadratic. Furthermore, since the SB method is just a special case of a two-split ADMM algorithm, it seems likely that the ADMM algorithm will be faster than the SB method if the AL penalty parameters are selected appropriately. To verify our analysis, we conduct a convergence rate analysis of the ADMM algorithm with two split variables for image restoration problems with quadratic data-fitting term and regularization term. Our analysis shows that the two-split ADMM algorithm can be faster than the SB method if the AL penalty parameter of the SB method is suboptimal.

The paper is organized as follows. In Section \ref{sec:jour-13-acp:sb_admm}, we prove the convergence of the SB method that allows inexact updates for regularized least-squares problems by showing the equivalence of the SB method and a convergent ADMM algorithm. To have a more concrete example and mathematically tractble analysis, Section \ref{sec:jour-13-acp:freq_analysis} shows a convergence rate analysis of an ADMM algorithm with two split variables for \emph{quadratically} regularized least-squares problems. Based on our convergence rate analysis, a discussion about parameter selection of ADMM algorithms in practical situations are shown in Section \ref{sec:jour-13-acp:para_selection}. Section \ref{sec:jour-13-acp:numerical_expt} demonstrates the experimental results supporting our analysis. Finally, we draw our conclusions in Section \ref{sec:jour-13-acp:conclusions}.

\section{The split Bregman method as an ADMM algorithm} \label{sec:jour-13-acp:sb_admm}
To show the convergence of the inexact SB method, we first consider another constrained minimization problem that is also equivalent to \eqref{eq:reprot-13-asc:ls_general_l1} but uses two split variables:
\begin{equation} \label{eq:jour-13-acp:eq_ls_general_l1_two_splits}
	\left(\hat{\mb{x}},\hat{\mb{u}},\hat{\mb{v}}\right)
	=
	\argmin
	{\mb{x},\mb{u},\mb{v}}
	{\ts\frac{1}{2}\norm{\mb{y}-\mb{u}}{2}^2+\fx{\Phi}{\mb{v}}}
	\text{ s.t. }
	\mb{u}=\mb{Ax},\mb{v}=\msb{\Theta}\mb{x} \, .
\end{equation}
The ADMM algorithm for this constrained minimization problem is \cite{afonso:11:aal}:
\begin{equation} \label{eq:jour-13-acp:admm_iterates}
	\begin{cases}
	\iter{\mb{x}}{k+1}
	=
	\argmin
	{\mb{x}}
	{\ts\frac{\rho}{2}\norm{\mb{Ax}-\iter{\mb{u}}{k}-\iter{\mb{d}}{k}}{2}^2+\ts\frac{\eta}{2}\norm{\msb{\Theta}\mb{x}-\iter{\mb{v}}{k}-\iter{\mb{e}}{k}}{2}^2} \\
	\iter{\mb{u}}{k+1}
	=
	\argmin
	{\mb{u}}
	{\ts\frac{1}{2}\norm{\mb{y}-\mb{u}}{2}^2+\ts\frac{\rho}{2}\norm{\mb{A}\iter{\mb{x}}{k+1}-\mb{u}-\iter{\mb{d}}{k}}{2}^2} \\
	\iter{\mb{v}}{k+1}
	=
	\argmin
	{\mb{v}}
	{\fx{\Phi}{\mb{v}}+\ts\frac{\eta}{2}\norm{\msb{\Theta}\iter{\mb{x}}{k+1}-\mb{v}-\iter{\mb{e}}{k}}{2}^2} \\
	\iter{\mb{d}}{k+1}
	=
	\iter{\mb{d}}{k}-\mb{A}\iter{\mb{x}}{k+1}+\iter{\mb{u}}{k+1} \\
	\iter{\mb{e}}{k+1}
	=
	\iter{\mb{e}}{k}-\msb{\Theta}\iter{\mb{x}}{k+1}+\iter{\mb{v}}{k+1} \, ,
	\end{cases}
\end{equation}
where $\mb{d}$ and $\mb{e}$ are the scaled Lagrange multipliers (i.e., dual variables) of the split variables $\mb{u}$ and $\mb{v}$, respectively, and $\rho>0$ and $\eta>0$ are the corresponding AL penalty parameters. By stacking $\mb{u}$ and $\mb{v}$, we can represent the equality constraint in \eqref{eq:jour-13-acp:eq_ls_general_l1_two_splits} more compactly as
\begin{equation} \label{eq:jour-13-acp:variable_splitting}
	\begin{bmatrix}
	\mb{u} \\ \mb{v}
	\end{bmatrix}
	=
	\underbrace{
	\begin{bmatrix}
	\mb{A} \\ \msb{\Theta}
	\end{bmatrix}}_{\mb{S}}
	\mb{x} \, .
\end{equation}
When $\mb{S}$ has full column rank, this ADMM algorithm \eqref{eq:jour-13-acp:admm_iterates} is convergent, even with inexact updates, providing the error in the inexact updates satisfies the conditions of \cite[Theorem 8]{eckstein:92:otd}. In many applications such as image restoration and X-ray CT image reconstruction, $\mb{A}'\mb{A}$ is a low-pass filter (but not necessarily shift-invariant). When $\msb{\Theta}=\mb{C}$ is the finite difference matrix, $\msb{\Theta}'\msb{\Theta}$ is the Laplacian, which is a high-pass filter. The non-zero vectors in the null space of $\msb{\Theta}'\msb{\Theta}$ are usually not in the null space of $\mb{A}'\mb{A}$, and vice versa, so the null space of $\mb{S}'\mb{S}=\mb{A}'\mb{A}+\msb{\Theta}'\msb{\Theta}$ is usually $\left\{\mb{0}\right\}$. That is, $\mb{S}$ usually has full column rank in applications like image restoration and X-ray CT image reconstruction! Therefore, \eqref{eq:jour-13-acp:admm_iterates} is a convergent ADMM algorithm that allows inexact updates for image restoration and X-ray CT image reconstruction according to \cite[Theorem 8]{eckstein:92:otd}. More specifically, $\iter{\mb{x}}{k+1}$ in \eqref{eq:jour-13-acp:admm_iterates} converges to $\hat{\mb{x}}$ if the error of the inner minimization problem (i.e., the $\ell_2$ distance between the iterate and the optimum of the inner problem) is absolutely summable.

Now, let's take a closer look at \eqref{eq:jour-13-acp:admm_iterates}. The $\mb{u}$-update in \eqref{eq:jour-13-acp:admm_iterates} has a closed-form solution
\begin{equation} \label{eq:jour-13-acp:closed_form_u_update}
	\iter{\mb{u}}{k+1}
	=
	\ts\frac{\rho}{\rho+1}\left(\mb{A}\iter{\mb{x}}{k+1}-\iter{\mb{d}}{k}\right)
	+
	\ts\frac{1}{\rho+1}\mb{y} \, .
\end{equation}
Combining with the $\mb{d}$-update in \eqref{eq:jour-13-acp:admm_iterates}, we have the identity
\begin{equation} \label{eq:jour-13-acp:u_d_identity}
	\iter{\mb{u}}{k+1}+\rho\iter{\mb{d}}{k+1}=\mb{y}
\end{equation}
if we initialize $\mb{d}$ as $\iter{\mb{d}}{0}=\rho^{-1}\left(\mb{y}-\iter{\mb{u}}{0}\right)$. Substituting \eqref{eq:jour-13-acp:u_d_identity} into \eqref{eq:jour-13-acp:admm_iterates}, we have the \emph{simplified} ADMM iterates:
\begin{equation} \label{eq:jour-13-acp:simplified_admm_iterates}
	\begin{cases}
	\iter{\mb{x}}{k+1}
	=
	\argmin
	{\mb{x}}
	{
	\begin{aligned}
	\ts\frac{\rho}{2}\norm{\mb{Ax}-\rho^{-1}\mb{y}-\left(1-\rho^{-1}\right)\iter{\mb{u}}{k}}{2}^2 \qquad \\
	+\ts\frac{\eta}{2}\norm{\msb{\Theta}\mb{x}-\iter{\mb{v}}{k}-\iter{\mb{e}}{k}}{2}^2
	\end{aligned}
	} \\
	\iter{\mb{u}}{k+1}
	=
	\ts\frac{\rho}{\rho+1}\mb{A}\iter{\mb{x}}{k+1}
	+
	\ts\frac{1}{\rho+1}\iter{\mb{u}}{k} \\
	\iter{\mb{v}}{k+1}
	=
	\argmin
	{\mb{v}}
	{\fx{\Phi}{\mb{v}}+\ts\frac{\eta}{2}\norm{\msb{\Theta}\iter{\mb{x}}{k+1}-\mb{v}-\iter{\mb{e}}{k}}{2}^2} \\
	\iter{\mb{e}}{k+1}
	=
	\iter{\mb{e}}{k}-\msb{\Theta}\iter{\mb{x}}{k+1}+\iter{\mb{v}}{k+1} \, .
	\end{cases}
\end{equation}
By comparing the SB method \eqref{eq:jour-13-acp:sb_iterates} and the simplified ADMM algorithm \eqref{eq:jour-13-acp:simplified_admm_iterates} side by side, we can easily find that they have common $\mb{v}$- and $\mb{e}$-updates. The $\mb{u}$-update in \eqref{eq:jour-13-acp:simplified_admm_iterates} can be seen as a perturbation of its $\mb{x}$-update. In fact, when $\rho=1$, the $\mb{x}$-update in \eqref{eq:jour-13-acp:simplified_admm_iterates} is independent of $\iter{\mb{u}}{k}$, and the simplified ADMM algorithm \eqref{eq:jour-13-acp:simplified_admm_iterates} reduces to the SB method \eqref{eq:jour-13-acp:sb_iterates}. In other words, \emph{the SB method is a convergent ADMM algorithm when we solve a regularized least-squares problem}, and this proves the convergence of the \emph{inexact} SB method for image restoration and X-ray CT image reconstruction provided $\mb{S}$ in \eqref{eq:jour-13-acp:variable_splitting} has full column rank, and the inner minimization error is absolute summable! Note that $\mb{S}$ has full column rank in many applications whereas $\msb{\Theta}$ often does not. This is the main difference between our new convergence condition and the conventional one.

\section{Convergence rate analysis of ADMM algorithms: the quadratic case} \label{sec:jour-13-acp:freq_analysis}
In the previous section, we showed that when the data-fitting term is quadratic, the SB method is a convergent ADMM algorithm, and therefore proved the convergence of the inexact SB method. Although the convergence of the SB method for general convex data-fitting term is still an open problem, the convergence proof in Section \ref{sec:jour-13-acp:sb_admm} is applicable to many popular image reconstruction problems. Note that the equivalence of the SB method and the ADMM algorithm holds for the choice $\rho=1$; however, the ADMM algorithm \eqref{eq:jour-13-acp:simplified_admm_iterates} is convergent for any $\rho>0$. Thus, it seems likely that the ADMM algorithm will be faster than the SB method if $\rho$ is selected appropriately. To have a more concrete example and mathematically tractble analysis, we analyze the convergence rate properties of \eqref{eq:jour-13-acp:simplified_admm_iterates} for a \emph{quadratically} regularized image restoration problem:
\begin{equation} \label{eq:jour-13-acp:ir_quad_reg}
	\hat{\mb{x}}
	=
	\argmin{\mb{x}}{\ts\frac{1}{2}\norm{\mb{y}-\mb{Ax}}{2}^2+\frac{\alpha}{2}\norm{\mb{Cx}}{2}^2} \, ,
\end{equation}
where $\mb{y}$ denotes the noisy blurred measurement of an image $\mb{x}$, degraded by a degradation matrix $\mb{A}$, $\alpha>0$ is the regularization parameter, and $\mb{C}$ denotes the tall \emph{masked} finite difference matrix in multiple directions. To simplify our analysis, we will further assume that both $\mb{A}'\mb{A}$ and $\mb{C}'\mb{C}$ are approximately block circulant with circulant blocks (BCCB), i.e., $\mb{A}'\mb{A}\approx\mb{U}\msb{\Lambda}\mb{U}'$ and $\mb{C}'\mb{C}\approx\mb{U}\msb{\Omega}\mb{U}'$, where $\msb{\Lambda}\teq\Diag{\lambda_i\geq0}$, $\msb{\Omega}\teq\Diag{\omega_i\geq0}$, and $\mb{U}$ denotes the normalized $2$D inverse DFT matrix.

Clearly, the quadratically regularized image restoration problem \eqref{eq:jour-13-acp:ir_quad_reg} is simply an instance of the regularized least-squares problem \eqref{eq:reprot-13-asc:ls_general_l1} with $\msb{\Theta}\teq\mb{C}$ and $\fx{\Phi}{\cdot}\teq\frac{\alpha}{2}\norm{\cdot}{2}^2$. Therefore, the simplified ADMM algorithm solving \eqref{eq:jour-13-acp:ir_quad_reg} is:
\begin{equation} \label{eq:jour-13-acp:admm_ir_quad_reg_u_v}
	\begin{cases}
	\iter{\mb{x}}{k+1}
	=
	\argmin
	{\mb{x}}
	{
	\begin{aligned}
	\ts\frac{\rho}{2}\norm{\mb{Ax}-\rho^{-1}\mb{y}-\left(1-\rho^{-1}\right)\iter{\mb{u}}{k}}{2}^2 \qquad \\
	+\ts\frac{\eta}{2}\norm{\mb{C}\mb{x}-\iter{\mb{v}}{k}-\iter{\mb{e}}{k}}{2}^2
	\end{aligned}
	} \\
	\iter{\mb{u}}{k+1}
	=
	\ts\frac{\rho}{\rho+1}\mb{A}\iter{\mb{x}}{k+1}
	+
	\ts\frac{1}{\rho+1}\iter{\mb{u}}{k} \\
	\iter{\mb{v}}{k+1}
	=
	\argmin{\mb{v}}{\ts\frac{\alpha}{2}\norm{\mb{v}}{2}^2+\frac{\eta}{2}\norm{\mb{C}\iter{\mb{x}}{k+1}-\mb{v}-\iter{\mb{e}}{k}}{2}^2} \\
	\iter{\mb{e}}{k+1}=\iter{\mb{e}}{k}-\mb{C}\iter{\mb{x}}{k+1}+\iter{\mb{v}}{k+1} \, .
	\end{cases}
\end{equation}
Furthermore, since $\Phi$ is quadratic, it has a linear proximal mapping, and therefore, the $\mb{v}$-update in \eqref{eq:jour-13-acp:admm_ir_quad_reg_u_v} has a closed-form solution
\begin{equation} \label{eq:jour-13-acp:closed_form_v_update}
	\iter{\mb{v}}{k+1}
	=
	\ts\frac{\eta}{\eta+\alpha}\left(\mb{C}\iter{\mb{x}}{k+1}-\iter{\mb{e}}{k}\right) \, .
\end{equation}
Again, using the same trick as before, we find that the dual variable $\mb{e}$ is also redundant, yielding the identity
\begin{equation} \label{eq:jour-13-acp:v_e_identity}
	\alpha\iter{\mb{v}}{k+1}+\eta\iter{\mb{e}}{k+1}=\mb{0}
\end{equation}
if we initialize $\mb{e}$ as $\iter{\mb{e}}{0}=-\alpha\eta^{-1}\iter{\mb{v}}{0}$. Substituting \eqref{eq:jour-13-acp:v_e_identity} into \eqref{eq:jour-13-acp:admm_ir_quad_reg_u_v}, the ADMM iterates \eqref{eq:jour-13-acp:admm_ir_quad_reg_u_v} simplify to:
\begin{equation} \label{eq:jour-13-acp:simplified_admm_ir_quad_reg_u_v}
	\begin{cases}
	\iter{\mb{x}}{k+1}
	=
	\left(\rho\mb{A}'\mb{A}+\eta\mb{C}'\mb{C}\right)^{-1}\big(\mb{A}'\mb{y}+(\rho-1)\mb{A}'\iter{\mb{u}}{k}+(\eta-\alpha)\mb{C}'\iter{\mb{v}}{k}\big) \\
	\iter{\mb{u}}{k+1}=\frac{\rho}{\rho+1}\mb{A}\iter{\mb{x}}{k+1}+\frac{1}{\rho+1}\iter{\mb{u}}{k} \\
	\iter{\mb{v}}{k+1}=\frac{\eta}{\eta+\alpha}\mb{C}\iter{\mb{x}}{k+1}+\frac{\alpha}{\eta+\alpha}\iter{\mb{v}}{k} \, .
	\end{cases}
\end{equation}
To further simplify \eqref{eq:jour-13-acp:simplified_admm_ir_quad_reg_u_v}, let's denote
\begin{equation} \label{eq:jour-13-acp:admm_substitutions}
	\begin{cases}
	\mb{s}\teq\left(\rho\mb{A}'\mb{A}+\eta\mb{C}'\mb{C}\right)^{-1}\mb{A}'\mb{y} \\
	\mb{P}\teq(\rho-1)\left(\rho\mb{A}'\mb{A}+\eta\mb{C}'\mb{C}\right)^{-1}\mb{A}' \\
	\mb{Q}\teq(\eta-\alpha)\left(\rho\mb{A}'\mb{A}+\eta\mb{C}'\mb{C}\right)^{-1}\mb{C}' \, .
	\end{cases}
\end{equation}
It follows that
\begin{equation} \label{eq:jour-13-acp:simplified_admm_ir_quad_reg_u_v_2}
	\begin{cases}
	\iter{\mb{x}}{k+1}=\mb{s}+\mb{P}\iter{\mb{u}}{k}+\mb{Q}\iter{\mb{v}}{k} \\
	\iter{\mb{u}}{k+1}=\frac{\rho}{\rho+1}\mb{A}\big(\mb{s}+\mb{P}\iter{\mb{u}}{k}+\mb{Q}\iter{\mb{v}}{k}\big)+\frac{1}{\rho+1}\iter{\mb{u}}{k} \\
	\iter{\mb{v}}{k+1}=\frac{\eta}{\eta+\alpha}\mb{C}\big(\mb{s}+\mb{P}\iter{\mb{u}}{k}+\mb{Q}\iter{\mb{v}}{k}\big)+\frac{\alpha}{\eta+\alpha}\iter{\mb{v}}{k} \, ,
	\end{cases}
\end{equation}
and we have the transition equation of the split variables:
\begin{equation} \label{eq:jour-13-acp:transition_mat_u_v}
	\left[
	\begin{gathered}
	\iter{\mb{u}}{k+1} \\ \iter{\mb{v}}{k+1}
	\end{gathered}
	\,\right]
	=
	\underbrace{
	\left[
	\begin{aligned}
	\begin{gathered}
	\ts\frac{\rho}{\rho+1}\mb{AP}+\frac{1}{\rho+1}\mb{I}_n \\
	\ts\frac{\eta}{\eta+\alpha}\mb{CP}
	\end{gathered}
	&&
	\begin{gathered}
	\ts\frac{\rho}{\rho+1}\mb{AQ} \\
	\ts\frac{\eta}{\eta+\alpha}\mb{CQ}+\frac{\alpha}{\eta+\alpha}\mb{I}_m
	\end{gathered}
	\end{aligned}
	\,\right]
	}_{\mb{G}}
	\left[
	\begin{gathered}
	\iter{\mb{u}}{k} \\ \iter{\mb{v}}{k}
	\end{gathered}
	\,\right]
	+
	\left[
	\begin{gathered}
	\ts\frac{\rho}{\rho+1}\mb{As} \\ \ts\frac{\eta}{\eta+\alpha}\mb{Cs}
	\end{gathered}
	\,\right] \, .
\end{equation}
Since we already know that the two-split ADMM algorithm \eqref{eq:jour-13-acp:admm_ir_quad_reg_u_v} is convergent if $\mb{A}$ and $\mb{C}$ have disjoint null spaces (except for $\left\{\mb{0}\right\}$) \cite[Theorem 8]{eckstein:92:otd}, the split variables in \eqref{eq:jour-13-acp:admm_ir_quad_reg_u_v} should converge linearly with rate of convergence $\fx{\varrho}{\mb{G}}$ \cite{ortega:iso:70,kelley:imf:95}, where $\fx{\varrho}{\cdot}$ denotes the spectral radius of a matrix. However, what we really care about is the convergence rate of $\mb{x}$. To find the convergence rate of $\mb{x}$, consider
\begin{align} \label{eq:jour-13-acp:primal_dual_eqn}
	& \,\,\,\,\,\,\,\,
	\iter{\mb{x}}{k+1}-\mb{s} \nonumber \\
	&=
	\left[
	\begin{aligned}
	\mb{P} && \mb{Q}
	\end{aligned}
	\,\right]
	\left[
	\begin{gathered}
	\iter{\mb{u}}{k} \\ \iter{\mb{v}}{k}
	\end{gathered}
	\,\right] \nonumber \\
	&=
	\left[
	\begin{aligned}
	\mb{P} && \mb{Q}
	\end{aligned}
	\,\right]
	\left(
	\left[
	\begin{aligned}
	\begin{gathered}
	\ts\frac{\rho}{\rho+1}\mb{AP}+\frac{1}{\rho+1}\mb{I}_n \\
	\ts\frac{\eta}{\eta+\alpha}\mb{CP}
	\end{gathered}
	&&
	\begin{gathered}
	\ts\frac{\rho}{\rho+1}\mb{AQ} \\
	\ts\frac{\eta}{\eta+\alpha}\mb{CQ}+\frac{\alpha}{\eta+\alpha}\mb{I}_m
	\end{gathered}
	\end{aligned}
	\,\right]
	\left[
	\begin{gathered}
	\iter{\mb{u}}{k-1} \\ \iter{\mb{v}}{k-1}
	\end{gathered}
	\,\right]
	+
	\left[
	\begin{gathered}
	\ts\frac{\rho}{\rho+1}\mb{As} \\ \ts\frac{\eta}{\eta+\alpha}\mb{Cs}
	\end{gathered}
	\,\right]
	\right) \nonumber \\
	&=
	\left(\ts\frac{\rho}{\rho+1}\mb{PA}+\frac{1}{\rho+1}\mb{I}_n+\frac{\eta}{\eta+\alpha}\mb{QC}\right)\mb{P}\iter{\mb{u}}{k-1}
	+
	\left(\ts\frac{\rho}{\rho+1}\mb{PA}+\frac{\eta}{\eta+\alpha}\mb{QC}+\frac{\alpha}{\eta+\alpha}\mb{I}_n\right)\mb{Q}\iter{\mb{v}}{k-1} \nonumber \\
	& \,\,\,\,\,\,~~~~~~~~~~~~~~~~~~~~~~~~~~~~~~~~~~~~~~~~~~~~~~~~~~~~~~~~~~~~~~~~~~~~~~~~~
	+
	\left(\ts\frac{\rho}{\rho+1}\mb{PA}+\frac{\eta}{\eta+\alpha}\mb{QC}\right)\mb{s} \nonumber \, .
\end{align}
Unfortunately, this is not a transition equation of $\mb{x}$ (or $\mb{x}-\mb{s}$), so we cannot find the linear convergence rate of $\mb{x}$ in general, except for three cases: (1) $\rho=1$, (2) $\eta=\alpha$, and (3) $\rho=\eta/\alpha$.

\subsection{Case I} \label{subsec:jour-13-acp:case_1}
When $\rho=1$, $\mb{P}$ becomes a zero matrix, and the split variable $\mb{u}$ is redundant. In fact, the two-split ADMM algorithm \eqref{eq:jour-13-acp:admm_ir_quad_reg_u_v} reduces to the SB method when $\rho=1$ as shown in Section \ref{sec:jour-13-acp:sb_admm}. In this case, we have
\begin{equation} \label{eq:jour-13-acp:transition_mat_v_sb}
	\iter{\mb{v}}{k+1}=\left(\ts\frac{\eta}{\eta+\alpha}\mb{CQ}+\frac{\alpha}{\eta+\alpha}\mb{I}_m\right)\iter{\mb{v}}{k}+\ts\frac{\eta}{\eta+\alpha}\mb{Cs}
\end{equation}
and
\begin{align} \label{eq:jour-13-acp:transition_mat_x_sb}
	& \,\,\,\,\,\,\,\,
	\iter{\mb{x}}{k+1}-\mb{s} \nonumber \\
	&=
	\left(\ts\frac{\eta}{\eta+\alpha}\mb{QC}+\frac{\alpha}{\eta+\alpha}\mb{I}_n\right)\mb{Q}\iter{\mb{v}}{k-1}+\ts\frac{\eta}{\eta+\alpha}\mb{QCs} \nonumber \\
	&=
	\underbrace{
	\left(
	\ts\frac{\eta}{\eta+\alpha}\mb{QC}+\frac{\alpha}{\eta+\alpha}\mb{I}_n
	\right)
	}_{\mb{H}_1}
	\big(\iter{\mb{x}}{k}-\mb{s}\big)+\ts\frac{\eta}{\eta+\alpha}\mb{QCs} \, .
\end{align}
Because the two-split ADMM algorithm \eqref{eq:jour-13-acp:admm_ir_quad_reg_u_v} is convergent, it follows that $\mb{x}$ converges linearly to the solution $\hat{\mb{x}}$ with rate $\fx{\varrho}{\mb{H}_1}$. Now, applying our BCCB approximations of $\mb{A}'\mb{A}$ and $\mb{C}'\mb{C}$, we can approximate the transition matrix $\mb{H}_1$ as
\begin{align} \label{eq:jour-13-acp:approx_H1}
	\mb{H}_1
	&=
	\ts\frac{\eta}{\eta+\alpha}\left((\eta-\alpha)\left(\mb{A}'\mb{A}+\eta\mb{C}'\mb{C}\right)^{-1}\mb{C}'\right)\mb{C}+\frac{\alpha}{\eta+\alpha}\mb{I}_n \nonumber \\
	&\approx
	\mb{U}\,
	\Diag{\ts\frac{\eta}{\eta+\alpha}\frac{(\eta-\alpha)\omega_i}{\lambda_i+\eta\omega_i}+\frac{\alpha}{\eta+\alpha}}
	\mb{U}' \nonumber \\
	&=
	\mb{U}\,
	\Diag{\ts\frac{\eta}{\eta+\alpha}\frac{\alpha\lambda_i+\eta^2\omega_i}{\eta\lambda_i+\eta^2\omega_i}}
	\mb{U}' \nonumber \\
	&=
	\mb{U}\,
	\Diag{\fx{s_1}{\delta_i}\teq\ts\frac{\eta}{\eta+\alpha}\frac{\alpha+\eta^2\delta_i}{\eta+\eta^2\delta_i}}
	\mb{U}' \, ,
\end{align}
where $\delta_i\teq\omega_i/\lambda_i\geq0$ is the ratio of the spectra of $\mb{C}'\mb{C}$ and $\mb{A}'\mb{A}$. Note that for any non-negative $\delta$, $\fx{s_1'}{\delta}=\frac{\eta}{\eta+\alpha}\frac{\eta-\alpha}{(\eta\delta+1)^2}$ is greater zero if $\eta>\alpha$, and it is less than zero if $\eta<\alpha$. When $\eta>\alpha$, $\fx{\varrho}{\mb{H}_1}$, i.e., $\nbvarmax{i}{\fx{s_1}{\delta_i}}$, is determined by $\delta_{\text{max}}$, and we can find the optimal AL penalty parameter
\begin{equation} \label{eq:jour-13-acp:opt_eta_sb_case_1}
	\eta^{\star}
	=
	\argmin{\eta}{\ts\frac{\eta}{\eta+\alpha}\frac{\alpha+\eta^2\delta_{\text{max}}}{\eta+\eta^2\delta_{\text{max}}}}
	=
	\sqrt{\alpha/\delta_{\text{max}}} \, .
\end{equation}
Note that \eqref{eq:jour-13-acp:opt_eta_sb_case_1} holds only if $\eta>\alpha$; therefore, $\eta^{\star}=\sqrt{\alpha/\delta_{\text{max}}}$ only if $\delta_{\text{max}}<\alpha^{-1}$. Similarly, when $\eta<\alpha$, $\fx{\varrho}{\mb{H}_1}$ is determined by $\delta_{\text{min}}$. Follow the same procedure, we have $\eta^{\star}=\sqrt{\alpha/\delta_{\text{min}}}$ only if $\delta_{\text{min}}>\alpha^{-1}$. Finally, for the case that $\delta_{\text{min}}<\alpha^{-1}<\delta_{\text{max}}$, $\eta^{\star}=\alpha$ because both $\eta^{\star}>\alpha$ and $\eta^{\star}<\alpha$ lead to a contradiction to the condition $\delta_{\text{min}}<\alpha^{-1}<\delta_{\text{max}}$. Summarizing, the optimal AL penalty parameter of the SB method for solving the quadratically regularized image restoration problem \eqref{eq:jour-13-acp:ir_quad_reg} is
\begin{equation} \label{eq:jour-13-acp:opt_eta_sb}
	\eta^{\star}=\ts\sqrt{\alpha/\gamma} \, ,
\end{equation}
where $\gamma\teq\text{median}\!\left\{\delta_{\text{min}},\delta_{\text{max}},\alpha^{-1}\right\}$.

\subsection{Case II} \label{subsec:jour-13-acp:case_2}
When $\eta=\alpha$, $\mb{Q}$ becomes a zero matrix, and the split variable $\mb{v}$ is redundant. Furthermore, we can easily prove that when $\eta=\alpha$, the two-split ADMM algorithm \eqref{eq:jour-13-acp:admm_ir_quad_reg_u_v} reduces to the alternating direction AL method that solves the constrained minimization problem:
\begin{equation} \label{eq:jour-13-acp:ir_quad_reg_u}
	\left(\hat{\mb{x}},\hat{\mb{u}}\right)
	=
	\argmin{\mb{x},\mb{u}}{\ts\frac{1}{2}\norm{\mb{y}-\mb{u}}{2}^2+\frac{\alpha}{2}\norm{\mb{Cx}}{2}^2}
	\text{ s.t. }
	\mb{u}=\mb{Ax}
\end{equation}
that is also equivalent to \eqref{eq:jour-13-acp:ir_quad_reg}. In this case, we have
\begin{equation} \label{eq:jour-13-acp:transition_mat_u_al}
	\iter{\mb{u}}{k+1}=\left(\ts\frac{\rho}{\rho+1}\mb{AP}+\frac{1}{\rho+1}\mb{I}_n\right)\iter{\mb{u}}{k}+\ts\frac{\rho}{\rho+1}\mb{As}
\end{equation}
and
\begin{align} \label{eq:jour-13-acp:transition_mat_x_al}
	& \,\,\,\,\,\,\,\,
	\iter{\mb{x}}{k+1}-\mb{s} \nonumber \\
	&=
	\left(\ts\frac{\rho}{\rho+1}\mb{PA}+\frac{1}{\rho+1}\mb{I}_n\right)\mb{P}\iter{\mb{u}}{k-1}+\ts\frac{\rho}{\rho+1}\mb{PAs} \nonumber \\
	&=
	\underbrace{
	\left(
	\ts\frac{\rho}{\rho+1}\mb{PA}+\frac{1}{\rho+1}\mb{I}_n
	\right)
	}_{\mb{H}_2}
	\big(\iter{\mb{x}}{k}-\mb{s}\big)+\ts\frac{\rho}{\rho+1}\mb{PAs} \, .
\end{align}
Follow the same trick, we can approximate the transition matrix $\mb{H}_2$ as
\begin{align} \label{eq:jour-13-acp:approx_H2}
	\mb{H}_2
	&=
	\ts\frac{\rho}{\rho+1}\left((\rho-1)\left(\rho\mb{A}'\mb{A}+\alpha\mb{C}'\mb{C}\right)^{-1}\mb{A}'\right)\mb{A}+\frac{1}{\rho+1}\mb{I}_n \nonumber \\
	&\approx
	\mb{U}\,
	\Diag{\ts\frac{\rho}{\rho+1}\frac{(\rho-1)\lambda_i}{\rho\lambda_i+\alpha\omega_i}+\frac{1}{\rho+1}}
	\mb{U}' \nonumber \\
	&=
	\mb{U}\,
	\Diag{\ts\frac{\rho}{\rho+1}\frac{\rho^2\lambda_i+\alpha\omega_i}{\rho^2\lambda_i+\alpha\rho\omega_i}}
	\mb{U}' \nonumber \\
	&=
	\mb{U}\,
	\Diag{\fx{s_2}{\delta_i}\teq\ts\frac{\rho}{\rho+1}\frac{\rho^2+\alpha\delta_i}{\rho^2+\alpha\rho\delta_i}}
	\mb{U}' \, ,
\end{align}
and the optimal AL penalty parameter $\rho^{\star}$ will be
\begin{equation} \label{eq:reprot-13-faa:opt_rho_al}
	\rho^{\star}=\ts\sqrt{\alpha\gamma} \, .
\end{equation}

\subsection{Case III} \label{subsec:jour-13-acp:case_3}
Finally, when $\rho=\eta/\alpha$, we have the identity $\frac{1}{\rho+1}=\frac{\alpha}{\eta+\alpha}$ and
\begin{align} \label{eq:jour-13-acp:transition_mat_x_admm}
	& \,\,\,\,\,\,\,\,
	\iter{\mb{x}}{k+1}-\mb{s} \nonumber \\
	&=
	\left(\ts\frac{\eta}{\eta+\alpha}\mb{PA}+\frac{\eta}{\eta+\alpha}\mb{QC}+\frac{\alpha}{\eta+\alpha}\mb{I}_n\right)\left(\mb{P}\iter{\mb{u}}{k-1}+\mb{Q}\iter{\mb{v}}{k-1}\right) \nonumber \\
	& \,\,\,\,\,\,~~~~~~~~~~~~~~~~~~~~~~~~~~~~~~~~~~~~~~~~~~~~~~~~~~~~~~\,+\ts\frac{\eta}{\eta+\alpha}\left(\mb{PA}+\mb{QC}\right)\mb{s} \nonumber \\
	&=
	\underbrace{
	\left(\ts\frac{\eta}{\eta+\alpha}\mb{PA}+\frac{\eta}{\eta+\alpha}\mb{QC}+\frac{\alpha}{\eta+\alpha}\mb{I}_n\right)
	}_{\mb{H}_3}
	\big(\iter{\mb{x}}{k}-\mb{s}\big)+\ts\frac{\eta}{\eta+\alpha}\left(\mb{PA}+\mb{QC}\right)\mb{s} \, .
\end{align}
The transition matrix $\mb{H}_3$ is approximately
\begin{align} \label{eq:jour-13-acp:approx_H3}
	\mb{H}_3
	&=
	\ts\frac{\eta}{\eta+\alpha}\left(\frac{\eta}{\alpha}\mb{A}'\mb{A}+\eta\mb{C}'\mb{C}\right)^{-1}\left((\frac{\eta}{\alpha}-1)\mb{A}'\mb{A}+(\eta-\alpha)\mb{C}'\mb{C}\right)+\frac{\alpha}{\eta+\alpha}\mb{I}_n \nonumber \\
	&\approx
	\mb{U}\,
	\Diag{\ts\frac{\eta}{\eta+\alpha}\left(\frac{(\eta-\alpha)\lambda_i+\alpha(\eta-\alpha)\omega_i}{\eta\lambda_i+\alpha\eta\omega_i}\right)+\frac{\alpha}{\eta+\alpha}}
	\mb{U}' \nonumber \\
	&=
	\mb{U}\,
	\Diag{\ts\frac{\eta}{\eta+\alpha}}
	\mb{U}' \nonumber \\
	&=
	\mb{U}\,
	\Diag{\fx{s_3}{\delta_i}\teq\ts\frac{\eta}{\eta+\alpha}}
	\mb{U}' \, .
\end{align}
Surprisingly, $\mb{H}_3$ has a uniform sprectrum, and $\fx{\varrho}{\mb{H}_3}=\eta/(\eta+\alpha)$. Theoretically, we can achieve arbitrarily fast asymptotic convergence rate in this quadratic case by choosing
\begin{equation} \label{eq:jour-13-acp:opt_eta_rho_admm}
	\eta^{\star}\approx0 \, .
\end{equation}
However, a smaller AL penalty parameter leads to a larger step size. When $\eta$ is too small, we might encounter overshoots at the beginning and oscillation as the algorithm proceeds. Therefore, in practice, $\eta^{\star}$ cannot be arbitrarily small.

\section{Parameter selection of ADMM algorithms for image restoration problems: the quadratic case} \label{sec:jour-13-acp:para_selection}
This section considers parameter selection of ADMM algorithms for image restoration in practical situations, where $\mb{A}'\mb{A}$ is a non-invertible low-pass filter, $\mb{C}'\mb{C}$ is a non-invertible high-pass filter, some frequency band is non-zero only for $\mb{A}'\mb{A}$ (such as the DC component), and some frequency band is non-zero only for $\mb{C}'\mb{C}$ (such as the extremely high frequency component). In this case, $\delta_i$ has an extremely huge dynamic range, i.e., $\delta_{\text{min}}\approx0$ and $\delta_{\text{max}}\approx\infty$. Therefore, for most cases, the optimal AL penalty parameter $\eta^{\star}$ of the SB method \eqref{eq:jour-13-acp:opt_eta_sb} will be $\alpha$. Furthermore, the optimal AL penalty parameter $\rho^{\star}$ of the two-split ADMM algorithm when $\eta=\alpha$ (Case II) is one, which then reverts to Case I, i.e., the SB method. Hence, once $\eta$ is chosen to be the optimal $\eta$ of the SB method, the optimal two-split ADMM algorithm is the SB method itself!

Now, consider the case that $\eta$ is suboptimal, i.e., $\eta\neq\alpha$. When $\eta>\alpha$, $\delta_{\text{max}}\approx\infty$ determines the asymptotic convergence rate. In this case,
\begin{equation} \label{eq:jour-13-acp:eta_greater_than_alpha}
	\fx{\varrho}{\mb{H}_1}\approx\ts\frac{\eta}{\eta+\alpha}=\fx{\varrho}{\mb{H}_3} \, ,
\end{equation}
which means that the SB method is no better than the two-split ADMM algorithm \eqref{eq:jour-13-acp:admm_ir_quad_reg_u_v} with $\rho=\eta/\alpha$ when $\eta$ is over-estimated. However, in practice, the SB method appears to converge a little bit faster because most frequency components in the SB method have convergence rate less than $\fx{\varrho}{\mb{H}_1}\approx\fx{\varrho}{\mb{H}_3}$. On the other hand, when $\eta<\alpha$, $\delta_{\text{min}}\approx0$ determines the convergence rate. In this case,
\begin{equation} \label{eq:jour-13-acp:eta_less_than_alpha}
	\fx{\varrho}{\mb{H}_1}\approx\ts\frac{\alpha}{\eta+\alpha}>\ts\frac{\eta}{\eta+\alpha}=\fx{\varrho}{\mb{H}_3} \, ,
\end{equation}
which means that the SB method is slower than the two-split ADMM algorithm \eqref{eq:jour-13-acp:admm_ir_quad_reg_u_v} with $\rho=\eta/\alpha$ when $\eta$ is under-estimated. In sum, $\fx{\varrho}{\mb{H}_1}\gtrsim\fx{\varrho}{\mb{H}_3}$ for any $\eta$. That is, the two-split ADMM algorithm \eqref{eq:jour-13-acp:admm_ir_quad_reg_u_v} with $\rho=\eta/\alpha$ is less sensitive to the choice of $\eta$ due to the additional split and converges faster than the SB method especially for small $\eta$ in most cases.

This analysis of the two-split ADMM algorithm might seem to be useless because we assume that we can solve the inner least-squares problem exactly and efficiently in our analysis, while the minimization problem \eqref{eq:jour-13-acp:ir_quad_reg} itself is a least-squares problem. In fact, if we initialize $\mb{d}$ and $\mb{e}$ in \eqref{eq:jour-13-acp:admm_ir_quad_reg_u_v} properly as mentioned before, the two-split ADMM algorithm should solve the minimization problem in \emph{one} iteration if we set $(\rho,\eta)=(1,\alpha)$, which happens to be the optimal SB method, as in \eqref{eq:jour-13-acp:simplified_admm_ir_quad_reg_u_v}. This does not contradict the non-zero ($1/2$) asymptotic convergence rate we showed in \eqref{eq:jour-13-acp:approx_H1}, \eqref{eq:jour-13-acp:approx_H2}, and \eqref{eq:jour-13-acp:approx_H3} because the $\mb{x}$-update just solves the original minimization problem fortuitously. The other split variables still follow the asymptotic convergence rate we derived before. The goal of this analysis was to show that ADMM algorithms can sometimes converge faster than the SB method, and the simple analysis might give some intuition about the parameter tuning for practical problems.

\section{Numerical experiments} \label{sec:jour-13-acp:numerical_expt}
\begin{figure*}
	\centering
	\includegraphics[width=0.8\textwidth]{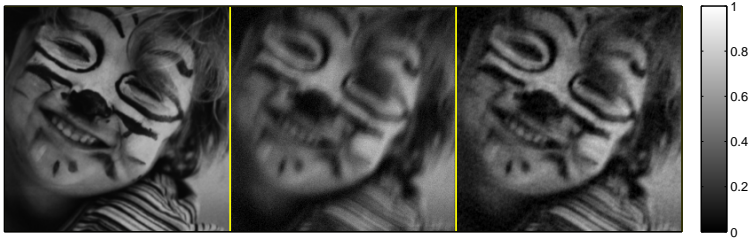}
	\caption{An image restoration problem instance: the true image (left), the noisy blurred image (middle), and the converged reference reconstruction (right).}
	\label{fig:jour-13-acp:ir_instance}
\end{figure*}

\begin{figure}
	\centering
	\subfigure[]{
	\includegraphics[width=0.75\textwidth]{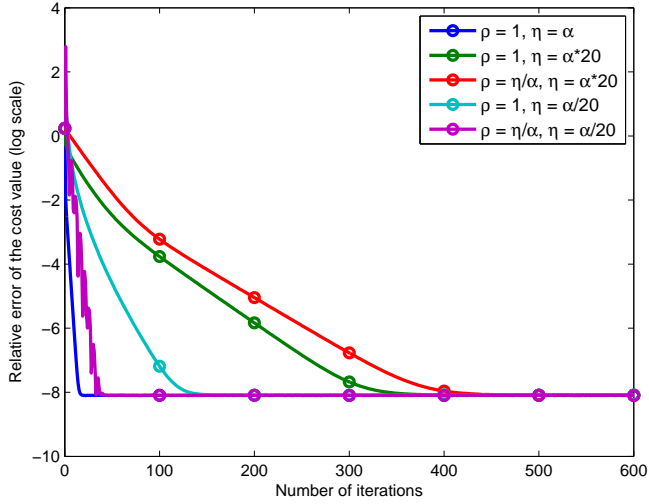}
	\label{fig:jour-13-acp:rel_error}
	}
	\subfigure[]{
	\includegraphics[width=0.75\textwidth]{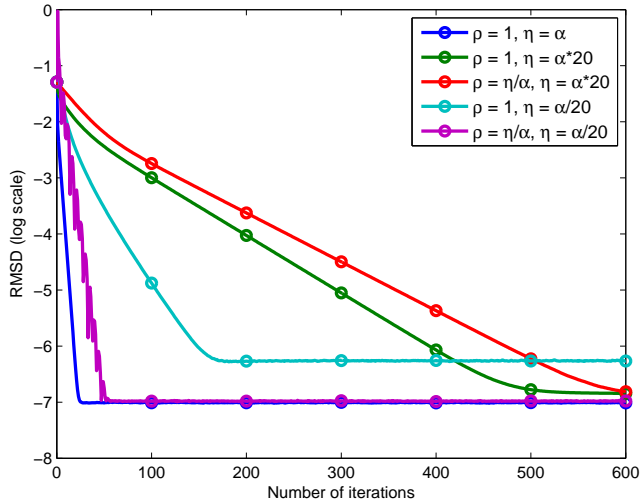}
	\label{fig:jour-13-acp:rms_diff}
	}
	\caption{The convergence rate curves with different parameter settings: (a) the relative error of cost value and (b) the RMS difference between the reconstructed image and the reference reconstruction as a function of the number of iterations.}
	\label{fig:jour-13-acp:conv_rate}
\end{figure}

In this section, we verify the convergence rate result and parameter selection discussed in the previous section using an image restoration problem with a quadratic regularizer. Figure \ref{fig:jour-13-acp:ir_instance} shows an image restoration problem instance: the true image (left), the noisy blurred image (middle), and the converged reference reconstruction (right). We use a quadratic roughness penalty as the regularizer where the regularization parameter $\alpha$ is choosen to be $2^{-4}$ for better noise-resolution tradeoff. Note that since a masked finite difference matrix (in horizontal and vertical directions) is used, we cannot solve the $\mb{x}$-update in \eqref{eq:jour-13-acp:admm_ir_quad_reg_u_v} efficiently using FFT. Instead, we solve it using PCG with an appropriate circulant preconditioner for three iterations. The inexact updates might affect the convergence rate but not very significantly thanks to the circulant preconditioner. Figure \ref{fig:jour-13-acp:conv_rate} shows the convergence rate curves (the relative error of cost value and RMS difference) of the two-split ADMM algorithm \eqref{eq:jour-13-acp:admm_ir_quad_reg_u_v} with different parameter settings. As can be seen in Figure \ref{fig:jour-13-acp:conv_rate}, all reconstructed images with different parameter settings converge to the solution with minimum cost value (up to the machine epsilon of the single-precision floating-point arithmetic). When $(\alpha,\eta)=(1,\alpha)$, the two-split ADMM algorithm, i.e., the optimal SB method, achieves the fastest convergence rate with no ripple. As mentioned before, with a proper initialization, the two-split ADMM algorithm with this parameter setting should converge immediately; the non-zero convergence rate comes from the inexact updates. When $\eta$ is over-estimated ($\eta=\alpha\!\times\!20$), the SB method ($\rho=1$) and the two-split ADMM algorithm with $\rho=\eta/\alpha=20$ exhibit similar slow convergence rate. When $\eta$ is under-estimated ($\eta=\alpha/20$), the SB method ($\rho=1$) is much slower than the two-split ADMM algorithm with $\rho=\eta/\alpha=1/20$. One might expect the ADMM algorithm with these parameters to converge with the same asymptotic convergence rate as the fastest two-split ADMM algorithm because $\eta$ is very small, but in fact it suffers from strong overshoots and oscillation due to the large step size as mentioned in Section \ref{subsec:jour-13-acp:case_3}.

\section{Conclusions} \label{sec:jour-13-acp:conclusions}
In this paper, we showed that for regularized least-squares problems, the split Bregman (SB) method is a convergent alternating direction method of multipliers (ADMM) for a broad class of regularizers. Therefore, the SB method has all the nice convergence properties of ADMM, such as the unconditional convergence with any augmented Lagrangian (AL) penalty parameters and inexact updates. Although the convergence of the SB method for general convex data-fitting terms is still an open problem, the proof in the paper is applicable to the most popular image reconstruction problems. To have a deeper understanding of the SB method and ADMM algorithms, we analyzed the convergence rate of the ADMM algorithm with two split variables for image restoration problems with a quadratic regularizer. According to our analysis, ADMM algorithms can sometimes converge faster than the SB method especially when the AL penalty parameter of the SB method is under-estimated. Although our analysis cannot be applied to image restoration problems with a non-quadratic edge-preserving regularizer, it gives insight on how to tune the AL penalty parameters for those pixels in which the cost function is almost quadratic, e.g., pixels in flat regions. As future works, we are interested in the convergence rate analysis of the SB method and ADMM algorithms with inexact updates, which might let us know how an inexact least-squares problem solver or an approximate proximal mapping of a more complicated proximal operator would affect the rate of convergence of these popular algorithms.

\section*{Acknowledgements}
This work was supported in part by NIH grant R01-HL-098686 and by an equipment donation from Intel.

\bibliographystyle{siam}
\bibliography{master}

\end{document}

%% file: def_siam.tex
\newcommand{\mb}[1]{\mathbf{#1}}
\newcommand{\msb}[1]{\boldsymbol{#1}}

\newcommand{\ts}{\textstyle}

\newcommand{\teq}{\triangleq}
\newcommand{\fx}[2]{#1\!\left(#2\right)}

\newcommand{\argmin}[2]{\underset{#1}{\text{argmin}}\left\{#2\right\}}

\newcommand{\nbvarmax}[2]{\underset{#1}{\text{max}}\,#2}

\newcommand{\norm}[2]{\left\|#1\right\|_{#2}}

\newcommand{\iter}[2]{#1^{\left(#2\right)}}

\newcommand{\Diag}[1]{\text{diag}\!\left\{#1\right\}}
